\input amstex
\input amsppt.sty
\magnification=\magstep1
\hsize=30truecc
\vsize=22.2truecm
\baselineskip=16truept
\NoBlackBoxes
\nologo
\pageno=1
\TagsOnRight
\topmatter
\def\Z{\Bbb Z}
\def\N{\Bbb N}

\def\l{\left}
\def\r{\right}
\def\bg{\bigg}
\def\({\bg(}
\def\[{\bg[}
\def\){\bg)}
\def\]{\bg]}
\def\t{\text}
\def\f{\frac}
\def\mo{\roman{mod}}

\def\em{\emptyset}
\def\se {\subseteq}

\def\sm{\setminus}

\def\eq{\equiv}
\def\ls{\leqslant}
\def\gs{\geqslant}

\def\Da{\Delta}
\def\da{\delta}

\def\Proof{\noindent{\it Proof}}

\def\Ack{\medskip\noindent {\bf Acknowledgments}}
\topmatter \hbox{J. Combin. Theory Ser. A, 115(2008), no.\,2,
345--353.}
\bigskip
\title {Determination of the two-color Rado number
for $a_1x_1+\cdots+a_{m}x_{m}=x_0$}\endtitle
\author {Song Guo$^1$ and Zhi-Wei Sun$^{2,*}$}\endauthor
\leftheadtext{Song Guo and Zhi-Wei Sun}
\rightheadtext{The two-color Rado number for $a_1x_1+\cdots+a_{m}x_{m}=x_0$}
\affil $^1$Department
of Mathematics, Huaiyin Teachers college
\\Huaian 223001, People's Republic of China
\\ {\tt guosong77\@sohu.com}\
\medskip
$^2$Department of Mathematics, Nanjing University
\\Nanjing 210093, People's Republic of China
\\ {\tt zwsun\@nju.edu.cn}
\\ http://math.nju.edu.cn/$^\sim$zwsun
\endaffil

\abstract For positive integers $a_1,a_2,\ldots,a_m$, we determine
the least positive integer $R(a_1,\ldots,a_m)$ such that for every
2-coloring of the set $[1,n]=\{1,\ldots,n\}$ with $n\gs
R(a_1,\ldots,a_m)$ there exists a monochromatic solution to the
equation $a_1x_1+\cdots+a_mx_m=x_0$ with $x_0,\ldots,x_m\in[1,n]$.
The precise value of $R(a_1,\ldots,a_m)$ is shown to be
$av^2+v-a$, where $a=\min\{a_1,\ldots,a_m\}$ and
$v=\sum_{i=1}^{m}a_i$. This confirms a conjecture of B. Hopkins
and D. Schaal.
\endabstract
\thanks  2000 {\it Mathematics Subject Classification}.
Primary 05D10; Secondary 11B75, 11D04.
\newline\indent *This author is responsible for communications,
and supported by the National Science Fund for
Distinguished Young Scholars (No. 10425103) and a Key Program of
NSF (No. 10331020) in China.
\endthanks
\endtopmatter
\document

\heading 1. Introduction\endheading

 Let $\N=\{0,1,2,\ldots\}$, and
$[a,b]=\{x\in \N : a\ls x\ls b\}$ for $a,b\in\N$.
For $k,n\in\Z^+=\{1,2,3,\ldots\}$, we call a function
 $\Da :[1,n]\to [0,k-1]$  a {\it $k$-coloring} of the set
 $[1,n]$, and $\Delta(i)$ the {\it color} of $i\in[1,n]$.
Given a $k$-coloring of the set $[1,n]$,
a solution to the linear diophantine equation
$$a_0x_0+a_1x_1+\cdots+a_mx_m=0\quad\ (a_0,a_1,\ldots,a_m\in\Z)$$
with $x_0,x_1,\ldots,x_m\in[1,n]$
is called {\it monochromatic} if
$\Da(x_0)=\Da(x_1)=\cdots=\Da(x_m)$.

Let $k\in\Z^+$. In 1916, I. Schur [S] proved that
if $n\in\Z^+$ is sufficiently large then for every {\it $k$-coloring} of
the set $[1,n]$, there exists a monochromatic solution to
$$x_1+x_2=x_0$$
with $x_0,x_1,x_2\in[1,n]$.

Let $k\in\Z^+$ and $a_0,a_1,\ldots,a_m\in\Z\sm\{0\}$. Provided that
$\sum_{i\in I}a_i=0$ for some $\em\not=I\se\{0,1,\ldots,m\}$,
R. Rado showed that for sufficiently large $n\in\Z^+$
the equation $a_0x_0+a_1x_1+\cdots+a_mx_m=0$
always has a monochromatic solution
when a $k$-coloring of $[1,n]$ is given;
the least value of such an $n$ is called the {\it $k$-color Rado number}
for the equation.
Since $-1+1=0$, Schur's theorem is a particular case of
Rado's result. The reader may consult the book [LR] by B. M. Landman and A. Robertson
for a survey of results on Rado numbers.

In this paper, we are interested in precise values of 2-color Rado numbers.
By a theorem of Rado [R], if $a_0,a_1,\ldots,a_m\in\Z$
contain both positive and negative integers and at least three of them are
nonzero, then the homogeneous linear equation
$$a_0x_0+a_1x_1+\cdots+a_mx_m=0$$
has a monochromatic solution with $x_0,\ldots,x_m\in[1,n]$
for any sufficiently large $n\in\Z^+$ and a $2$-coloring of $[1,n]$.
In particular, if $a_1,\ldots,a_{m}\in\Z^+\ (m\gs2)$ then
there is a least positive integer $n_0=R(a_1,\ldots,a_{m})$
such that for any $n\gs n_0$ and a $2$-coloring of $[1,n]$
the diophantine equation
$$a_1x_1+\cdots+a_{m}x_{m}=x_0\tag1.0$$
always has a monochromatic solution with $x_0,\ldots,x_m\in[1,n]$.

 In 1982, A. Beutelspacher and W. Brestovansky [BB]
proved that the 2-color Rado number $R(1,\ldots,1)$
for the equation $x_1+\cdots+x_{m}=x_0\ (m\gs2)$ is $m^2+m-1$.
In 1991, H. L. Abbott [A] extended this by showing that
for the equation $$a(x_1+\cdots+x_{m})=x_0\ \ \ \ (a\in\Z^+\ \t{and}\ m\gs2)$$
the corresponding 2-color Rado number $R(a,\ldots,a)$
is $a^3m^2+am-a$; that $R(a,\ldots,a)\gs a^3m^2+am-a$
was first obtained by L. Funar [F], who conjectured the equality.
In 2001, S. Jones and D. Schaal [JS] proved that
if $a_1,\ldots,a_{m}\in\Z^+\ (m\gs2)$ and
$\min\{a_1,\ldots,a_{m}\}=1$ then
$R(a_1,\ldots,a_{m})=b^2+3b+1$ where $b=a_1+\cdots+a_{m}-1$;
this result actually appeared earlier in Funar [F].

 In 2005 B. Hopkins and D. Schaal [HS] showed the following result.

 \proclaim {Theorem 1.0} Let $m\gs2$ be an integer and let $a_1,\ldots,a_{m}\in\Z^+$.
Then
$$R(a,b)\gs R(a_1,\ldots,a_{m})\gs a(a+b)^2+b, \tag1.1$$
where
$$a=\min\{a_1,\ldots,a_{m}\}\ \ \t{and}\ \ b=\sum_{i=1}^{m}a_i-a.\tag1.2$$
\endproclaim

Hopkins and Schaal ([HS]) conjectured further that the two inequalities in (1.1)
are actually equalities and verified this in the case $a=2$.

In this paper we confirm the conjecture of Hopkins and Schaal;
namely, we establish the following theorem.

\proclaim{Theorem 1.1} Let $m\gs2$ be an integer and let $a_1,\ldots,a_{m}\in\Z^+$.
Then
$$R(a_1,\ldots,a_{m})=a(a+b)^2+b, \tag1.3$$
where $a$ and $b$ are as in $(1.2)$.
\endproclaim

By Theorem 1.1, if $a_1,\ldots,a_{m}\in\Z^+$ and $n\gs av^2+v-a$
with $a=\min\{a_1,\ldots,a_{m}\}$ and $v=a_1+\cdots+a_{m}$, then
for any $X\se[1,n]$ either there are $x_1,\ldots,x_m\in X$
such that $\sum_{i=1}^{m}a_ix_i\in X$
or there are $x_1,\ldots,x_{m}\in [1,n]\sm X$
such that $\sum_{i=1}^{m}a_ix_i\in[1,n]\sm X$.
\smallskip

In the next section we reduce Theorem 1.1 to the following weaker version.
\proclaim{Theorem 1.2} Let $a,b,n\in\Z^+$, $a\ls b$ and
$n\gs av^2+b$ with $v=a+b$. Suppose that $b(b-1)\not\eq0\ (\mo\ a)$
and
 $\Da:[1,n]\to[0,1]$ is a $2$-coloring of $[1,n]$
with $\Da(1)=0$ and $\Da(a)=\Da(b)=\da\in[0,1]$.
Then there is a monochromatic solution
to the equation
$$ax+by=z\quad\ (x,y,z\in[1,n]).\tag1.4$$
\endproclaim

In Sections 3 and 4 we will prove Theorem 1.2
in the cases $\da=0$ and $\da=1$ respectively.

\heading {2. Reduction of Theorem 1.1 to Theorem 1.2}\endheading

Let us first give a key lemma which will be used in Sections 2--4.

\proclaim{Lemma 2.1} Let $k,l,n\in\Z^+$ with $l<n$, and let $\Da:[1,n]\to[0,1]$
be a $2$-coloring of $[1,n]$. Suppose that $kx+ly=z$ has no
monochromatic solution with $x,y,z\in[1,n]$. Assume also that
$u$ is an element of $[1,n-l]$ with $\Da(u)=\da$ and $\Da(u+l)=1-\da$.

{\rm (i)} If  $w\in\Z^+$, $w\ls(n-ku)/l$ and $\Da(w)=\da$,
then $\Da(w-hk)=\da$ whenever $h\in \N$ and $w-hk>0$.

{\rm (ii)} If $w\in[1,n]$ and $\Da(w)=1-\da$, then $\Da(w+hk)=1-\da$
whenever $h\in \N$ and $w+hk\ls(n-ku)/l$.

\endproclaim
\Proof. It suffices to handle the case $h=1$, since we can
consider $w\mp(h-1)k$ instead of $w$ if $h>1$.

(i) As $\Da(u)=\Da(w)=\da$ and $w\ls (n-ku)/l$, we have $\Da(ku+lw)=1-\da$.
By $\Da(u+l)=1-\da$ and
$k(u+l)+l(w-k)=ku+lw$, if $w-k>0$ then $\Da(w-k)=\da$.

(ii) Since $\Da(u+l)=\Da(w)=1-\da$ and $(w+k)l+ku\ls n$, we have $\Da(k(u+l)+lw)=\da$.
Note that $\Da(u)=\da$ and
$ku+l(w+k)=k(u+l)+lw$. So $\Da(w+k)=1-\da$.

 The proof of Lemma 2.1 is now complete. \qed
 \medskip

 Now we deduce Theorem 1.1 from Theorem 1.2.

\medskip
\noindent{\it Proof of Theorem 1.1}. By Theorem 1.0, it suffices to show that
$R(a,b)\ls av^2+b$,
where $v=a+b$. Since $m\gs 2$, we have $a\ls b$.

Let $n\gs av^2+b$ be an integer and let $\Da:[1,n]\to[0,1]$
be a 2-coloring of $[1,n]$.
Without loss of generality, we may assume that $\Da(1)=0$.
Suppose, for contradiction, that
there doesn't exist any
monochromatic solution to the equation (1.4).

Since $a\cdot1+b\cdot1=v$, we have $\Da(v)\not=\Da(1)=0$,
and hence
$$\Da(v)=1. \tag 2.1$$
Similarly, as
$av+bv=v^2$, we must have
$$\Da(v^2)=0\ \ \t{and}\ \ \Da(av^2+b\cdot1)=1. \tag 2.2$$

{\bf Claim 2.1}. $\Da(a)=\Da(b)\not=\Da(av)=\Da(bv)$.

As $aa+ba=av$ and $ab+bb=bv$, we have
$$\Da(av)\not=\Da(a)\ \ \t{and}\ \ \Da(bv)\not=\Da(b).$$
If $\Delta(a)\not=\Da(b)$, then
$$\Da(av)=\Da(b)\not=\Da(a)=\Da(bv)$$
and hence
$$\Da(a)=\Da(ab+b(av))=\Da(abv+ab)=\Da(a(bv)+ba)=\Da(b),$$
which contradicts $\Da(a)\not=\Da(b)$.
\medskip

Below, we let $\da=\Da(a)=\Da(b)$ and hence $\Da(av)=\Da(bv)=1-\da$.

In view of Claim 2.1 and
Theorem 1.2, $a$ divides $b(b-1)$ since (1.4) has no monochromatic
solution.
\medskip

{\bf Claim 2.2}. $\Da(ab+bv+(1-\da)av)=0$.

Recall that $\Da(v)=1$ by (2.1). If $\da=1$, then $\Da(b)=1=\Da(v)$,
and hence $\Da(ab+bv)=0$. When $\da=0$, we have
$\Da(b)=0<\Da(b+a)=\Da(v)=1$, and hence $\Da(v+b)=1$ by Lemma 2.1(ii)
(with $k=u=b,\,l=a$ and $w=v$) since $v+b=a+2b\ls(n-b^2)/a$;
therefore, $\Da(a(v+b)+bv)=0$. This completes the proof of Claim 2.2.
\medskip

Observe that
$$ab(v-1)+ab+b+(1-\da)av\ls abv+b+av\ls av^2+b\ls n.$$

{\bf Claim 2.3}. For every $i=1,\ldots,a$ we have
$$\Da(ib(v-1)+ab+b+(1-\da)av)=0.\tag2.3$$

When $i=1$, (2.3) holds by Claim 2.2. Now let $1<i\ls a$ and assume that
(2.3) holds with $i$ replaced by $i-1$.
Then
$$\align&\Da\l(a\l((i-1)\f{b(b-1)}a+ib+(1-\da)v\r)+b\cdot1\r)
\\=&\Da((i-1)b(v-1)+ab+b+(1-\da)av)=0=\Da(1)
\endalign$$
by the induction hypothesis.
Therefore,
$$\Da\l((i-1)\f{b(b-1)}a+ib+(1-\da)v\r)=1=\Da(v),$$
and hence
$$\align&\Da(ib(v-1)+ab+b+(1-\da)av)
\\=&\Da\l(a\l((i-1)\f{b(b-1)}a+ib+(1-\da)v\r)+bv\r)
=0.\endalign$$
This concludes the induction proof of Claim 2.3.
\medskip

Putting $i=a$ in (2.3) we find
that $$\Da\l(abv+b+(1-\da)av\r)=0=\Da(1).$$ If $\da=1$, then
$\Da(a(bv)+b\cdot1)=0=\Da(1)$, and hence $\Da(bv)=1=\Da(b)$, which
is impossible by Claim 2.1.
Thus $\da=0$ and
$$\Da(aa+b(av+a+1))=\Da(abv+b+av)=0=\Da(a).$$
It follows that $\Da(av+a+1)=1$. Also, if $a=1$ then
$\Da(av^2+b)=\Da(abv+b+av)=0$.
Since $\Da(av^2+b)=1$ by (2.2), and
$$a(av-b)+b(av+a+1)=av^2+b,$$
we must have $a\gs2$ and $\Da(av-b)=0$. As
$\Da(b)=0<\Da(b+a)=1$ and
$$av-b=v^2-b(v+1)<v^2-b(b-1)\ls v^2-\f {b(b-1)}a\ls\f{n-b^2}a,$$
we have $\Da(a^2+b)=\Da(av-b-(a-2)b)=0$ by Lemma 2.1(i)
with $k=u=b,\,l=a$ and $w=av-b$. However,
$\Da(a^2+b)=\Da(aa+b\cdot1)=1$ since $\Da(a)=0=\Da(1)$, so we get a
contradiction. This completes the proof. \qed

\heading{3. Proof of Theorem 1.2 with $\da=0$}\endheading

To prove Theorem 1.2 in the case $\da=0$,
we should deduce a contradiction under the assumption that
(1.4) has no monochromatic solution. Recall the condition $\Da(1)=\Da(a)=\Da(b)=0$.
It is clear that $\Da(a\cdot1+b\cdot1)\not=\Da(1)=0$.

Note that $a(v-1)+b(v-1)=v^2-v\ls av^2+b\ls n$. We make the
following claim first.

{\bf Claim 3.1}. $\Da(ai+bj)=1$ for any $i,j\in[1,a]$.

Since $\Da(a)=0<\Da(a+b)=1$ and
$$v+(i-1)a\ls a^2+b\ls \f{ab^2+2a^2b+a^3-a^2}b<\f{av^2+b-a^2}b\ls\f{n-a^2}b,$$
we have $\Da(ai+b)=1$ by Lemma 2.1(ii) with $k=u=a,\,l=b$ and
$w=v$. Similarly, as
$$(ai+b)+b(j-1)\ls a^2+ab=av=v^2-bv<v^2-b\f{b-1}a\ls\f{n-b^2}a,$$
by Lemma 2.1(ii) with $k=u=b,\,l=a$ and $w=ai+b$ we get that
$$\Da (ai+bj)=\Da(ai+b+b(j-1))=\Da(ai+b)=1. $$
This proves Claim 3.1.

{\bf Claim 3.2}. $\Da(c)=0$ for any $c\in[1,v-1]$.

Suppose that $c\in[b+1,v-1]$ and $\Da(c)=1$. Then
$\Da(av+bc)=0=\Da(a)$ since $\Da(v)=1=\Da(c)$.
Therefore,
$$\Da(a(av+bc)+ba)=1.$$
Clearly,
$$a(av+bc)+ba=a(a^2+b(c-b+1))+b(a^2+ba),$$
and $\Da(a^2+b(c-b+1))=1=\Da(a^2+ba)$ by Claim 3.1. Thus we get a
monochromatic solution to (1.4), contradicting our assumption.
So, $\Da(c)=0$ for all $c\in[b+1,v-1]$.

Now let $c\in[1,b]$. Then there is $\bar c\in[b,v-1]$ such that
$\bar c-c=ha$ for some $h\in\N$ (e.g., $\bar c=c$ when $c=b$). Recall that
$\Da(a)=0<\Da(a+b)=\Da(v)=1$ and also $\Da(b)=0$. As $\Da(\bar c)=0$ and
$$\bar c<v<\f{v^2-a}b<\f{av^2+b-a^2}b\ls\f{n-a^2}b,$$
by Lemma 2.1(i) with $k=u=a,\,l=b$ and $w=\bar c$, we have
$\Da(c)=\Da(\bar c-ha)=0$. This concludes the proof of Claim 3.2.

{\bf Claim 3.3}. $\Da(ai+bj)=1$ for any $i,j\in[1,v-1]$.

By Claim 3.2 we have $\Da(i)=\Da(j)=0$. Thus $\Da(ai+bj)=1$ since
(1.4) has no monochromatic solution. So Claim 3.3 holds.
\bigskip

Let $d$ be the greatest common divisor of $a$ and $b$.
Since $a\nmid b$, we have $d<a<b$, hence both $a'=a/d$ and $b'=b/d$
are greater than one.
By elementary number theory, there is $s\in[1,b'-1]$
such that $a's\eq 1\ (\mo\ b')$. Since $1<a's<a'b'$, we have
$t=(a's-1)/b'\in[1,a'-1]$ and $b't<a'b'\ls av$.
Observe that
$$a(av+b's)+b(av-b't)=av(a+b)+b'd=av^2+b\ls n.$$
As $\Da(v^2)=\Da(av+bv)\not=\Da(v)=1$, we have $\Da(v^2)=0=\Da(1)$ and hence
$\Da(av^2+b\cdot1)=1$. Therefore,
$$\Da(av+b's)=0\ \ \t{or}\ \ \Da(av-b't)=0.\tag3.1$$
Since $a+s,a-t\in[1,v-1]$, we have
$$\Da(av+bs)=\Da(a^2+b(a+s))=1=\Da(a^2+b(a-t))=\Da(av-bt)$$
by Claim 3.3, which contradicts (3.1) if $b=b'$. So $b'\not=b$, and hence $d>1$.

In view of (3.1), we distinguish two cases.

{\it Case} 3.1. $\Da(av+b's)=0$.

 Choose $s_1\in\Z^+$ such that $1\ls
as_1-b't\ls a$. Since $as_1\ls a+b't\ls a+b(a-1)\ls ab$, we have
$s_1\ls b$. Clearly, $\Da(aa+ba)=\Da(as_1+b\cdot1)=1$ by Claim
3.3, and
$$a(a^2+ab)+b(as_1+b)\ls a^2v+b(a+ba)\ls av^2+b\ls n.$$
Therefore,
$$\Da(a(a^2+ab)+b(as_1+b))=0.$$
However,
$$a(a^2+ab)+b(as_1+b)=a(av+b's)+b(as_1-b't+b-1)$$
and $\Da(as_1-b't+b-1)=0$ by Claim 3.2. This contradicts the
assumption that (1.4) has no monochromatic solution.

{\it Case} 3.2. $\Da(av-b't)=0$.

 Choose $s_2\in\Z$ so that $0\ls a't-as_2\ls a-1$.
Clearly, $0\ls s_2\ls t\ls a'-1<a-1$. With the help of Claim 3.2,
$\Da(a't-as_2+b)=0=\Da(av-b't)$. Since
$$a(av-b't)+b(a't-as_2+b)=a^2v-abs_2+b^2\ls av^2+b\ls n,$$
we have $\Da(a^2v-abs_2+b^2)=1$. Observe that
$$a^2v-abs_2+b^2=a(a^2+b)+b(a(a-1-s_2)+b)$$
and $\Da(a^2+b)=\Da(a(a-1-s_2)+b)=1$ by Claim 3.3. So we get a
monochromatic solution to (1.4), contradicting our assumption.

\heading{4. Proof of Theorem 1.2 with $\da=1$}\endheading

Assume the conditions of Theorem 1.2 with $\da=1$,
and that (1.4) doesn't have a monochromatic solution.
Our goal is to deduce a contradiction.

Since $\Da(a)=\Da(b)=\da=1$, $av=aa+ba$ and $bv=ab+bb$, we have
$$\Da(av)=\Da(bv)=0.\tag4.1$$
Thus there is a positive multiple $u_1\ls b(v-1)$ of $b$ such that
$\Da(u_1)=1$ and $\Da(u_1+b)=0$; also there is a positive multiple
$u_2\ls a(v-1)$ of $a$ such that $\Da(u_2)=1$ and $\Da(u_2+a)=0$.

Observe that
$$a^2+a+1<a^2\cdot\f vb+a+1=\f{(av^2+b)-ab(v-1)}{b}\ls\f{n-au_1}{b}.$$
As $\Da(1)=0$ and $1+a<a^2+a+1$, we have
$\Da(1+a)=0$ by Lemma 2.1(ii) with
$k=a,\,l=b,\,u=u_1$ and $w=1$.
Thus,
$$\Da(av+v)=\Da(a(a+1)+b(a+1))=1.\tag 4.2$$

{\bf Claim 4.1}. $\Da(a^2+a)=1\Rightarrow \Da(a)=\Da(2a)=\cdots=\Da(a^2)=1$.

Recall that $\Da(u_1)=1>\Da(u_1+b)=0$ and $a^2+a<(n-au_1)/b$.
 By Lemma 2.1(i) with
$k=a,\,l=b,\,u=u_1$ and $w=a^2+a$, if $\Da(a^2+a)=1$ then
$\Da(a^2+a-ha)=1$ for all $h=0,\ldots,a$. This proves Claim 4.1.

{\bf Claim 4.2}. For $w\in[1,n]$ and $h\in\N$ with $w+hb\ls
av+b$, we have $\Da(w)=0\Rightarrow\Da(w+hb)=0$.

 Note that
$$av+b<av+b+\f ba=\f{(av^2+b)-ab(v-1)}{a}\ls\f{n-bu_2}{a}.$$
So we get Claim 4.2 by applying Lemma 2.1(ii) with
$k=b,\,l=a$ and $u=u_2$.
\medskip

Write $b=aq+r$ with $q,r\in\N$ and $r<a$.
Since $a\ls b$ and $a\nmid b(b-1)$, we have $q\gs 1$ and $r\gs2$.

{\bf Claim 4.3}. $\Da(r)=0\Longrightarrow \Da(a^2)=0$.

Assume that $\Da(r)=0$.
As $\Da(r+aq)=\Da(b)=1$, there is $u_3\in\{r,r+a,\ldots,r+a(q-1)\}$
such that $\Da(u_3)=0$ and $\Da(u_3+a)=1$.
Since $\Da(av)=0$ (cf. (4.1)) and
$$av=v^2-bv<v^2-b^2<v^2+b-b(b-1)\ls\f{(av^2+b)-b(b-a)}a\ls\f{n-bu_3}a,$$
we have $\Da(a^2)=\Da(av-ab)=0$ by Lemma 2.1(i) with $k=b,\,l=a,\,u=u_3$ and $w=av$.

{\bf Claim 4.4}. $\Da(r)=\Da(ar)=1\Longrightarrow \Da(av+a)=0$.

Suppose that $\Da(r)=\Da(ar)=1$. Then $\Da(vr)=\Da(ar+br)=0$.
So there is $u_4\in\{ar, ar+b,\ldots,ar+(r-1)b\}$ such that
$\Da(u_4)=1$ and $\Da(u_4+b)=0$. Since
$\Da(av)=0$  by (4.1), and
$$av+a<av+a(a-r)\f vb=av\f{v-r}b<\f{(av^2+b)-a(vr-b)}b\ls\f{n-au_4}b,$$
we have $\Da(av+a)=0$ by Lemma 2.1(ii) with $k=a,\,l=b,\,u=u_4$ and $w=av$.

{\bf Claim 4.5}. $\Delta(av+a)=0$.

Clearly, $(a^2+a)+ab=av+a\ls av+b$. If $\Da(a^2+a)=0$, then we have
$\Da(av+a)=\Da((a^2+a)+ab)=0$ by applying Claim 4.2 with $w=a^2+a$ and $h=a$.
In the case $\Da(a^2+a)=1$, by Claim 4.1 we have $\Da(a^2)=1=\Da(ar)$, hence
$\Da(r)=\Da(ar)=1$ by Claim 4.3 and $\Da(av+a)=0$ by Claim 4.4.
\medskip

{\bf Claim 4.6}. There exists $u\in[1,ab-a]$ such that $\Da(u)=1$ and $\Da(u+a)=0$.

As $a$ does not divide $b$, the greatest common divisor $d$ of $a$ and $b$
is smaller than $a$, and hence $1<a'=a/d<b'=b/d$.
If $\Da(db)=0$, then we have
$$\Da(ab)=\Da(db+(a-d)b)=0<1=\Da(a)$$
by applying Claim 4.2 with $w=db$ and $h=a-d$,
hence there is $u\in\{a,2a,\ldots,(b-1)a\}$
such that $\Da(u)=1$ and $\Da(u+a)=0$.
Below we work under the condition $\Da(db)=1$.

{\it Case} 4.1. $\Da(d)=1$.

In this case, $d>1$ since $\Da(d)\not=\Da(1)$.
Note that $\Da(dv)=\Da(ad+bd)=1-\Da(d)=0$.
As $\Da(db)=1$,
for some $u\in\{db,db+a,\ldots, db+(d-1)a\}$ we have
$\Da(u)=1>\Da(u+a)=0$. Note that $a'b'-a'-b'=(a'-1)(b'-1)-1\gs0$ and
$$u\ls dv-a=d^2(a'+b')-a\ls d^2a'b'-a=ab-a.$$

{\it Case} 4.2. $\Da(d)=0$.

Choose $s\in[0,b-1]$ such that $as\eq d\ (\mo\ b)$. Clearly
$s\not=0,1$ since $d<a<b$. For $t=(as-d)/b$, we have $0<t<a$. As
$\Da(d)=0$ and $d+bt=as<ab\ls av+b$, we have $\Da(as)=\Da(d+bt)=0$
by Claim 4.2 with $w=d$ and $h=t$. Recall that $\Da(a)=1$. So
there is $u\in\{a,2a,\ldots,(s-1)a\}$ such that $\Da(u)=1$ and
$\Da(u+a)=0$. Clearly, $u\ls(s-1)a<ab-a$. This concludes the proof
of Claim 4.6.
\medskip

Let $u$ be as required in Claim 4.6.
Then
$$av+v=v^2-(b-1)v\ls v^2-b(b-1)<\f{(av^2+b)-b(ab-a)}a\ls\f{n-bu}a.$$
Recall that $\Da(av+v)=1$ by (4.2).
Thus $\Da(av+a)=\Da((av+v)-b)=1$ by Lemma 2.1(i) with $k=b,\,l=a,$ and $w=av+v$.
This contradicts Claim 4.5 and we are done.

\Ack. The authors are very grateful to Prof. B. Hopkins for
reading the initial version carefully, Prof. L. Funar for
informing the authors about the references [A] and [F], and the
two referees for many helpful comments.


\medskip

\widestnumber\key{BB}

\Refs

\ref\key A\by H. L. Abbott\paper On a conjecture of Funar concerning
generalized sum-free sets
\jour Nieuw Arch. Wisk (4)\vol 9\yr 1991\pages 249--252\endref

\ref\key BB\by A. Beutelspacher, W. Brestovansky
\paper Generalized Schur numbers\jour
in: Combinatorial Theory (Schloss Rauischholzhausen, 1982),
Lecture Notes in Math., 969, Springer, New York, 1982, pp. 30--38\endref

\ref\key F\by L. Funar\paper Generalized sum-free sets of integers
\jour Nieuw Arch. Wisk (4)\vol 8\yr 1990\pages 49--54\endref

\ref\key HS\by B. Hopkins and D. Schaal \paper On Rado numbers for
$\sum_{i=1}^{m-1}a_ix_i=x_m$
\jour Adv. in Appl. Math. \vol 35\yr 2005\pages 433--441\endref

\ref\key JS\by S. Jones and D. Schaal \paper Some 2-color Rado numbers
\jour Congr. Numer.\vol 152\yr 2001\pages 197--199\endref

\ref\key LR\by B. M. Landman and A. Robertson
\book Ramsey Theory on the Integers\publ American Mathematical Society,
Providence, RI, 2004\endref

\ref\key R\by R. Rado\paper Studien zur Kombinatorik\jour Math. Z.\vol
36\yr 1933 \pages 242--280\endref

\ref\key S\by I. Schur\paper \"Uber die Kongruenz $x^m+y^m\eq z^m
\ (\mo\ p)$ \jour Jahresber. Deutsch. Math.-Verein.\vol 25\yr 1916\pages
114--116\endref

\endRefs
\enddocument